\documentclass[letterpaper]{article}

\usepackage{amsmath,amsthm,bbm,amssymb,amsfonts,enumerate,latexsym,mathrsfs}

\RequirePackage[left=1in,right=1in,top=1in,bottom=1in]{geometry}

\newcommand{\G}{\mathbb{G}}

\newtheorem{theorem}{Theorem}[section]

\theoremstyle{definition}
\newtheorem{definition}[theorem]{Definition}
\newtheorem{remark}[theorem]{Remark}
\newtheorem{example}[theorem]{Example}
\newtheorem*{example*}{Example}

\DeclareMathOperator{\SL}{SL}

\DeclareMathOperator{\Hom}{Hom}

\DeclareMathOperator{\Lie}{Lie}
\DeclareMathOperator{\id}{id}
\DeclareMathOperator{\End}{End}

\DeclareMathOperator{\Spec}{Spec}

\DeclareMathOperator{\Ext}{Ext}

\DeclareMathOperator{\ad}{ad}

\begin{document}

\title{The big projective module as a nearby cycles sheaf}
\author{Justin Campbell}
\maketitle

\begin{abstract}

We give a new geometric construction of the big projective module in the principal block of the BGG category $\mathscr{O}$, or rather the corresponding $\mathscr{D}$-module on the flag variety. Namely, given a one-parameter family of nondegenerate additive characters of the unipotent radical of a Borel subgroup which degenerate to the trivial character, there is a corresponding one-parameter family of Whittaker sheaves. We show that the unipotent nearby cycles functor applied to this family yields the big projective $\mathscr{D}$-module.

\end{abstract}

\section{Introduction}

Let $k$ be an algebraically closed field of characteristic zero and $G$ a reductive group over $k$. Fix a Borel subgroup $B \subset G$ with unipotent radical $N \subset B$, and choose a splitting $T := B/N \to B$, so that we may speak of the opposite Borel $B^-$ and its unipotent radical $N^-$. Write $\Delta$ for the set of simple roots of $G$ with respect to $B$.

We study the abelian category $\mathscr{D}(G/B)^N$ of $N$-equivariant holonomic $\mathscr{D}$-modules on the flag variety $G/B$, which is equivalent to the principal block of the BGG category $\mathscr{O}$ under Beilinson-Bernstein localization (see e.g. \cite{BB} and \cite{S}). The category $\mathscr{D}(G/B)^N$ has finitely many simple objects, labeled by elements of the Weyl group $W$, which we will denote by $\mathscr{L}_w$ for $w \in W$. Write $C_w = NwB/B$ for the Schubert cell corresponding to $w$, and $X_w = \overline{C_w}$ for the Schubert variety.  Recall that $\mathscr{L}_w$ is the IC (intersection cohomology) sheaf on $X_w$, pushed forward to $G/B$. Let $\mathscr{M}_w$ and $\mathscr{M}_w^{\vee}$ be the $!$- and $*$-pushforwards to $G/B$, respectively, of the IC sheaf on $C_w$. Denote by $\mathscr{P}_w$ a projective cover of $\mathscr{L}_w$, i.e. an indecomposable projective object of $\mathscr{D}(G/B)^N$ which maps nontrivially to $\mathscr{L}_w$. Recall that $\mathscr{P}_w$ is unique up to non-canonical isomorphism.

Consider $\mathscr{L}_e$, which is the delta sheaf at the closed $N$-orbit in $G/B$. Its projective cover $\mathscr{P}_e$ is the longest indecomposable projective object of $\mathscr{D}(G/B)^N$; it is often referred to as the \emph{big projective}.

The sheaf $\mathscr{P}_e$ can be constructed from a Whittaker sheaf by averaging as follows. This construction appears in Lemma 4.4.11 of \cite{BY}, as well as Proposition 14.3.1 of \cite{FG} (see also the related \cite{BBM}). Fix a generic additive character $\psi^- : N^- \to \mathbb{G}_a$ and write $e^{\psi^-} := (\psi^-)^!e^x[1-\dim N^-]$, where $e^x$ is the exponential $\mathscr{D}$-module on $\mathbb{G}_a$. The \emph{Whittaker sheaf} $\mathscr{W}(\psi^-)$ is obtained by pushing $e^{\psi^-}$ forward along the canonical open embedding $N^- \to G/B$ (the $*$- and $!$-pushforwards are the same, i.e. the extension is clean). Then $\mathscr{P}_e$ is isomorphic to the sheaf obtained by averaging $\mathscr{W}(\psi^-)$ against the action of $N$. The $!$- and $*$-averaging functors give the same result thanks to the Verdier self-duality of $\mathscr{P}_e$, although the canonical morphism from the $!$-average to the $*$-average is not an isomorphism. 

Our construction also begins with a Whittaker sheaf and yields $\mathscr{P}_e$, but proceeds in a different fashion (the author does not know how to formally link the two constructions). Namely, we apply the geometric Jacquet functor of \cite{ENV}. This operation proceeds by choosing a dominant regular cocharacter of $T$ and a generic character $\psi : N \to \mathbb{G}_a$, which gives rise to a $\mathbb{G}_m$-equivariant one-parameter family of Whittaker sheaves on $G/B$. We will prove that the unipotent nearby cycles of such a family is isomorphic to $\mathscr{P}_e$.

One can interpret this theorem as follows. The category of nondegenerate Whittaker (meaning $e^{\psi}$-equivariant) sheaves on $G/B$ is equivalent to the category of finite-dimensional vector spaces, being generated by $\mathscr{W}(\psi)$. On the other hand, the most singular block of category $\mathscr{O}$ is also equivalent to finite-dimensional vector spaces, and translation to the principal block sends the generator of the singular block to the big projective. So our theorem says that the geometric Jacquet functor, applied to nondegenerate Whittaker sheaves on $G/B$, provides a geometric model for translation from the most singular block of category $\mathscr{O}$ to the principal block.

As explained in Remark 14 of \cite{W}, Theorem 4.4.1 of \cite{BY} implies that a category of degenerate Whittaker sheaves, where $\psi$ is not necessarily generic, is equivalent to a possibly singular block of category $\mathscr{O}$. In future work the author intends to show that translation out of the wall corresponds to a geometric Jacquet functor between the corresponding categories of degenerate Whittaker sheaves.

An analogous construction can be made on Drinfeld's compactification of the moduli stack of (twisted) $N$-bundles on a smooth projective curve. Namely, one applies the geometric Jacquet functor to an appropriately defined nondegenerate Whittaker sheaf. In his forthcoming Ph.D. thesis, the author will study the novel tilting sheaf on Drinfeld's compactification obtained by this procedure.

\emph{Acknowledgements.} First and foremost I thank my advisor Dennis Gaitsgory, without whom this paper would not exist. I am grateful to Sam Raskin for many helpful conversations, as well as feedback which substantially simplified the proof of the theorem. I also thank the anonymous referee who pointed out several errors in an earlier draft.

\emph{Notations and conventions.} For any variety $X$ over $k$ we write $\mathscr{D}(X)$ for the abelian category of holonomic $\mathscr{D}$-modules on $X$. If $f : X \to Y$ is a smooth morphism, we write \[ f^{\Delta} := f^![\dim Y - \dim X] : \mathscr{D}(Y) \to \mathscr{D}(X) \] for the cohomologically normalized pullback functor.

\section{The Whittaker sheaf}

\label{whit}

Let $H$ be an algebraic group over $k$ and write $\mu : H \times H \to H$ for the group operation. A \emph{character sheaf} (sometimes called a multiplicative local system) on $H$ is a line bundle with connection $\chi$ on $H$ that satisfies $\mu^{\Delta}\chi \cong \chi \boxtimes \chi$. If $\varphi  : K \to H$ is a homomorphism, then $\varphi^{\Delta}\chi$ is a character sheaf on $K$.

On $\mathbb{G}_a$ the basic character sheaf is the exponential $\mathscr{D}$-module $e^x$. We can pull it back along an additive character $\varphi : H \to \mathbb{G}_a$ to obtain a character sheaf $e^{\varphi} := \varphi^{\Delta}e^x$ on $H$.

If $H$ acts on a variety $X$ we denote by $\mathscr{D}(X)^{\chi}$ the abelian category of $\chi$-equivariant holonomic $\mathscr{D}$-modules on $X$; when $\chi$ is the trivial character sheaf this is just the category $\mathscr{D}(X)^H$ of $H$-equivariant holonomic $\mathscr{D}$-modules on $X$. An object of $\mathscr{D}(X)^{\chi}$ consists of $\mathscr{F} \in \mathscr{D}(X)$ together with an isomorphism $\chi \boxtimes \mathscr{F} \tilde{\to} a^{\Delta}\mathscr{F}$ satisfying the usual cocycle condition on $H \times H \times X$ (here $a : H \times X \to X$ is the action map).

All of the above also makes sense in families, i.e. over base schemes other than $\Spec k$. We will only need to work with families over $\mathbb{A}^1$. For example, if $p : \mathbb{G}_a \times \mathbb{A}^1 \to \mathbb{G}_a$ denotes the projection, then $p^{\Delta}e^x$ is a character sheaf on the constant group scheme $\mathbb{G}_a \times \mathbb{A}^1$.

Resume the notation of the introduction. We fix an additive character \[ \psi : N \longrightarrow \mathbb{G}_a \] which is \emph{nondegenerate} in the sense that the induced homomorphism \[ N/[N,N] = \bigoplus_{\alpha \in \Delta} N_{\alpha} \to \mathbb{G}_a \] is nontrivial on every simple root group $N_{\alpha}$.

Recall that the big cell $C_{w_0}$ is an $N$-torsor that is trivialized by our choice of maximal torus $T \to B$. Write $j : N \to G/B$ for the resulting inclusion. 

\begin{definition}

The \emph{Whittaker sheaf attached to} $\psi$ is defined by $\mathscr{W}(\psi) := j_*e^{\psi}$.

\end{definition}

The canonical morphism $j_!e^{\psi} \to \mathscr{W}(\psi)$ is an isomorphism, i.e. the extension is clean. Indeed, one checks that $\psi$ is nontrivial on the stabilizer of any point not in the big cell. Since $\mathscr{W}(\psi)$ is $e^{\psi}$-equivariant by construction, this implies that its $*$-restriction to $C_w$ vanishes for any $w \neq w_0$.

Now we extend $\mathscr{W}(\psi)$ to a one-parameter family of character sheaves using a choice of dominant regular coweight $\gamma : \mathbb{G}_m \to T$. The assumption of dominance implies that the $\gamma$-conjugation action $(s,n) \mapsto \ad_{\gamma(s)}(n)$ of $\mathbb{G}_m$ on $N$ extends uniquely to an action of the multiplicative monoid $\mathbb{A}^1$. Thus $\gamma$ induces a homomorphism of constant group schemes $\psi_{\gamma} : N \times \mathbb{A}^1 \to \mathbb{G}_a \times \mathbb{A}^1$ given on points by \[ \psi_{\gamma}(n,s) = (\psi(\ad_{\gamma(s)}(n)),s). \] Now pulling back the exponential gives the character sheaf $e^{\psi_{\gamma}} := \psi_{\gamma}^{\Delta}p^{\Delta}e^x$ on $N \times \mathbb{A}^1$, where as before $p : \mathbb{G}_a \times \mathbb{A}^1 \to \mathbb{G}_a$ is the projection. Since $\gamma$ was chosen regular the $!$-fiber over $0 \in \mathbb{A}^1$ of $e^{\psi_{\gamma}}$ is, up to shift, the trivial character sheaf on $N$.

Finally, we obtain the desired one-parameter family $\mathscr{W}(\psi,\gamma) := (j \times \id_{\mathbb{A}^1})_*e^{\psi_{\gamma}}$, an $e^{\psi_{\gamma}}$-equivariant $\mathscr{D}$-module on $(G/B) \times \mathbb{A}^1$. The extension is clean away from $0 \in \mathbb{A}^1$, but the $!$-restriction of $\mathscr{W}(\psi,\gamma)$ to the fiber over $0$ is (up to shift) $\mathscr{M}^{\vee}_{w_0}$.

\section{Nearby cycles}

\label{nc}

Let $X$ be a variety over $k$. Write $\mathring{\mathbb{A}}^1 := \mathbb{A}^1 \setminus \{ 0 \}$. Recall the existence of the \emph{unipotent nearby cycles functor} \[ \Psi : \mathscr{D}(X \times \mathring{\mathbb{A}}^1) \longrightarrow \mathscr{D}(X), \] which has the following properties:
\begin{enumerate}

\item $\Psi$ is exact,
\item $\Psi$ is compatible with restriction to open subvarieties and pushforward along proper morphisms,
\item if $\mathscr{F} \in \mathscr{D}(X \times \mathbb{A}^1)$ is lisse, then $\Psi(\mathscr{F}|_{X \times \mathring{\mathbb{A}}^1}) = i^!\mathscr{F}[1]$ where $i : X \to X \times \mathbb{A}^1$ is the inclusion of the fiber over $0$,
\item $\Psi$ naturally lifts to a functor taking values in the category whose objects are pairs $(\mathscr{G},m)$ where $\mathscr{G}$ is a $\mathscr{D}$-module on $X$ and $m$ is a nilpotent endomorphism of $\mathscr{G}$.

\end{enumerate}

For any $\mathscr{F} \in \mathscr{D}(X \times \mathring{\mathbb{A}}^1)$ the resulting nilpotent endomorphism of $\Psi(\mathscr{F})$ is called its \emph{monodromy}.

See \cite{B} for a construction of $\Psi$, which allows for nonconstant families. We will only need the case of a constant family with $X = G/B$ or a subvariety of $G/B$.

Now fix an algebraic group $H$ acting on $X$ and a character sheaf $\chi$ on the constant group scheme $H \times \mathbb{A}^1$. Write $\mathring{\chi}$ for the restriction of $\chi$ to $H \times \mathring{\mathbb{A}}^1$ and write $\chi_0 := \chi|^!_{H \times \{ 0 \}}[1]$. Then it follows from Beilinson's construction in \emph{loc. cit.} that $\Psi$ lifts to a functor of equivariant $\mathscr{D}$-modules \[ \Psi : \mathscr{D}(X \times \mathring{\mathbb{A}}^1)^{\mathring{\chi}} \longrightarrow \mathscr{D}(X)^{\chi_0}. \]

\section{The theorem}

In our situation, we can view unipotent nearby cycles as a functor \[ \Psi : \mathscr{D}((G/B) \times \mathbb{A}^1)^{e^{\psi_{\gamma}}} \longrightarrow \mathscr{D}(G/B)^N \] by first restricting to $(G/B) \times \mathring{\mathbb{A}}^1$, which we omit from the notation.

\begin{theorem}

We have $\Psi(\mathscr{W}(\psi,\gamma)) \cong \mathscr{P}_e$.

\label{theorem}

\end{theorem}

\begin{remark}

The theorem shows that the $\mathscr{D}$-module $\Psi(\mathscr{W}(\psi,\gamma))$ does not depend on $\gamma$, but in fact $\gamma$ determines the monodromy endomorphism. Namely, note that $\mathscr{P}_e$ is weakly $T$-equivariant and so $\mathfrak{h} := \Lie T$ acts on $\mathscr{P}_e$ by endomorphisms. Fixing an isomorphism $\varphi : \Psi(\mathscr{W}(\psi,\gamma)) \tilde{\to} \mathscr{P}_e$, the action of $-d\gamma \in \mathfrak{h}$ on $\mathscr{P}_e$ corresponds to the monodromy endomorphism of $\Psi(\mathscr{W}(\psi,\gamma))$ by Claim 2 in the appendix to \cite{AB} (this makes sense independently of $\varphi$ because $\End(\mathscr{P}_e)$ is commutative, see e.g. \cite{Ber}).

\end{remark}

\begin{example}

Consider the case $G = \SL_2$, so that $G/B = \mathbb{P}^1$. Since $N = \G_a$ we can take $\psi = \id$, so that $\mathscr{W}(\psi) = j_*e^x$ where $j : \mathbb{A}^1 \to \mathbb{P}^1$ is the inclusion. Moreover $T = \mathbb{G}_m$ and $\gamma$ is the $n^{\text{th}}$ power map on $\mathbb{G}_m$ for some $n > 0$.

Let $s = w_0 \in W = S_2$ be the nontrivial (and longest) element of the Weyl group. Then $\mathscr{L}_s $ is the IC sheaf on $\mathbb{P}^1$, $\mathscr{P}_s \cong \mathscr{M}_s$, and $\mathscr{L}_e$ is the delta sheaf at $\infty$. The sheaf $\mathscr{P}_e$ has a composition series \[ 0 = \mathscr{F}_0 \subset \mathscr{F}_1 \subset \mathscr{F}_2 \subset \mathscr{F}_3 = \mathscr{P}_e \] with $\mathscr{F}_1 \cong \mathscr{L}_e$, $\mathscr{F}_2/\mathscr{F}_1 \cong \mathscr{L}_s$, and $\mathscr{P}_e/\mathscr{F}_2 \cong \mathscr{L}_e$. The algebra of endomorphisms of $\mathscr{P}_e$ is canonically identified with $k[\epsilon]/(\epsilon^2)$, and the monodromy endomorphism of $\Psi(\mathscr{W}(\psi,\gamma))$ corresponds to $-n \epsilon$.

\label{sl2}

\end{example}

\begin{proof}[Proof of the theorem]

We will write $\Psi := \Psi(\mathscr{W}(\psi,\gamma))$ for brevity. The proof proceeds in two steps: we will first show that for all $w \neq e$ we have $\Ext^i(\Psi,\mathscr{L}_w) = 0$ for all $i$, then prove that $\Hom(\Psi,\mathscr{L}_e) = k$ and $\Ext^i(\Psi,\mathscr{L}_e) = 0$ for $i \neq 0$.

We claim that for any $w \neq e$ there exists a simple parabolic subgroup $P \subset G$ (in particular, $P \neq B$) such that $\mathscr{L}_w \cong \pi^{\Delta}\mathscr{F}$ for some object $\mathscr{F} \in 
\mathscr{D}(G/P)^N$, where $\pi : G/B \to G/P$ is the projection. Indeed, if $\alpha$ is a simple root such that $\ell(ws_{\alpha}) = \ell(w) - 1$, then $\overline{NwB}$ is stable under the right action of $P = B \sqcup Bs_{\alpha}B$. The claim follows because $\mathscr{L}_w$ is pushed forward from the IC sheaf on $X_w$.

In this paragraph we use the de Rham pushforward functor $\pi_* = \pi_!$, which takes values in the derived category of $\mathscr{D}$-modules on $G/P$ with holonomic cohomologies. We will show that $\pi_*\Psi = 0$, so that \[ \Ext^i(\Psi,\mathscr{L}_w) \cong \Ext^i(\pi_!\Psi,\mathscr{F}[-1]) = 0 \] for all $i$. Since $\Psi$ commutes with pushforward along proper morphisms, we need only prove that $\widetilde{\pi}_*\mathscr{W}(\psi,\gamma) = 0$, where $\widetilde{\pi} : (G/B) \times \mathring{\mathbb{A}}^1 \to (G/P) \times \mathring{\mathbb{A}}^1$ is the projection (again we omit the restriction to $(G/B) \times \mathring{\mathbb{A}}^1$ from our notation). Because $\mathscr{W}(\psi,\gamma)$ is cleanly extended from $C_{w_0} \times \mathring{\mathbb{A}}^1$, the sheaf $\widetilde{\pi}_*\mathscr{W}(\psi,\gamma)$ is $*$-extended from the $N \times \mathring{\mathbb{A}}^1$-orbit $(Nw_0P/P) \times \mathring{\mathbb{A}}^1$. Thus by the $e^{\psi_{\gamma}}$-equivariance it suffices to check that $i_{w_0}^!\widetilde{\pi}_*\mathscr{W}(\psi,\gamma) = 0$, where $i_{w_0} : \mathring{\mathbb{A}}^1 \to (G/P) \times \mathring{\mathbb{A}}^1$ is the constant section $w_0P$. One has $w_0Pw_0^{-1} \cap N = N_{-w_0\alpha}$, so the desired vanishing follows from base change. This is because for any $t \neq 0$, the character $\ad_{\gamma(t)}(\psi)$ is nontrivial when restricted to $N_{-w_0\alpha}$, and a nontrivial exponential $\mathscr{D}$-module has vanishing de Rham cohomology.

Finally, we show that $\Ext^i(\Psi,\mathscr{L}_e) = 0$ if $i \neq 0$ and $\Hom(\Psi,\mathscr{L}_e) = k$. Recall that there is a surjection $\mathscr{M}_{w_0}^{\vee} \to \mathscr{L}_e$ whose kernel does not have $\mathscr{L}_e$ as a subquotient, which implies that \[ \Ext^i(\Psi,\mathscr{M}_{w_0}^{\vee}) \tilde{\longrightarrow} \Ext^i(\Psi,\mathscr{L}_e) \] for all $i$. But $\mathscr{W}(\psi,\gamma)$ is lisse on $C_{w_0} \times \mathbb{A}^1$, so that properties 2 and 3 of $\Psi$ from Section \ref{nc} imply that $\Psi|_{C_{w_0}} = \mathscr{W}(\psi,\gamma)|^!_{C_{w_0} \times \{ 0 \}}[1]$. The latter is the IC (i.e. cohomologically normalized constant) sheaf on $C_{w_0}$, whence the claim.

\end{proof}


\begin{thebibliography}{10}

\bibitem{AB}
Arkhipov, S. and Bezrukavnikov, R.: Perverse sheaves on affine flags and Langlands dual group, Israel Journal of Mathematics \textbf{170}(1), 135-183 (2009)

\bibitem{B}
Beilinson, A.: How to glue perverse sheaves, in K-theory, arithmetic, and geometry, Lect. Notes Math., vol. 1289. Springer, Berlin (1987)

\bibitem{BB}
Beilinson, A. and Bernstein, J.: Localisation de $\mathfrak{g}$-modules, C. R. Acad. Sci. Paris, S\'{e}r. 1 \textbf{292}, 15-18 (1981)

\bibitem{Ber}
Bernstein, J. Trace in categories, in Operator algebras, Unitary representations, Enveloping algebras and Invariant theory, Progr. Math. \textbf{92}, 417-424. Birkh\"{a}user, Boston (1990)

\bibitem{BBM}
Bezrukavnikov, R., Braverman, A. and Mirkovic, I.: Some results about geometric Whittaker model, Adv. Math. \textbf{186}, 143-152 (2004)

\bibitem{BY}
Bezrukavnikov, R. and Yun, Z.: On Koszul duality for Kac-Moody groups, Represent. Theory \textbf{17}, 1-98 (2013)

\bibitem{ENV}
Emerton, M., Nadler, D., and Vilonen, K.: A geometric Jacquet functor, Duke Mathematical Journal \textbf{125}(2), 267-278 (2004)

\bibitem{FG}
Frenkel, E. and Gaitsgory, D.: Local geometric Langlands correspondence and affine Kac-Moody algebras, in Algebraic geometry and number theory. Springer, Berlin (2006)

\bibitem{S}
Soergel, W.: \emph{\'{E}quivalences de certaines classes de} $\mathfrak{g}$\emph{-modules}, C. R. Acad. Sci. Paris, S\'{e}r. 1 \textbf{303}, 725-728 (1986)

\bibitem{W}
Webster, B.: Singular blocks of parabolic category O and finite W-algebras, Journal of Pure and Applied Algebra \textbf{215}(12), 2797-2804 (2011)

\end{thebibliography}
\end{document}